\documentclass{article}

\usepackage[style=iso]{datetime2}

\usepackage{float}

\usepackage{enumitem} 
\usepackage{xspace} 
\usepackage{xcolor} 
\usepackage{algorithm,algpseudocode}

\usepackage{amsmath,amssymb}
\usepackage{amsthm}
\usepackage{bm}
\usepackage{pdfsync}
\usepackage{subfigure}
\usepackage{color}
\usepackage[english]{babel}

\usepackage[T1]{fontenc}

\setlist{itemsep=-4pt,topsep=0pt}
\algnewcommand{\algorithmicgoto}{\textbf{go to}}%
\algnewcommand{\Goto}[1]{\algorithmicgoto~step~\ref{#1}}%

\algdef{SE}[DOWHILE]{Do}{DoWhile}{\algorithmicdo}[1]{\algorithmicwhile\ #1}
\algnewcommand{\IIf}[1]{\State\algorithmicif\ #1\ \algorithmicthen}
\algnewcommand{\ElseIIf}[1]{\algorithmicelse\ #1}
\algnewcommand{\ElseI}[1]{\algorithmicelse\ #1}
\algnewcommand{\EndIIf}{\unskip\ \algorithmicend\ \algorithmicif}

\theoremstyle{plain}
\newtheorem{theorem}{Theorem}[section]

\theoremstyle{definition}
\newtheorem{example}[theorem]{Example}
\newtheorem{definition}[theorem]{Definition}
\newtheorem{remark}[theorem]{Remark}

\newcommand{\maxentry}[1]{|| #1 ||_{\max}}
\newcommand \Hrow {H_{\text{row}}}
\newcommand \Hcol {H_{\text{col}}}

\newcommand \ie {\textit{i.e.}}
\newcommand \eg {\textit{e.g.}}

\newcommand \valuation {\nu}  

\newcommand \Mat {\mathop{\rm Mat}}
\newcommand \adj {\mathop{\rm adj}}
\newcommand \softO {O^\sim}  

\newcommand \CC {{\mathbb C}}

\newcommand \NN {{\mathbb N}}

\newcommand \RR {{\mathbb R}}

\newcommand \ZZ {{\mathbb Z}}

\def\cocoa{\mbox{\rm
C\kern-.13em o\kern-.07 em C\kern-.13em o\kern-.15em A}}
\def\apcocoa{\mbox{\rm
A\kern-0.13em p\kern -0.07em C\kern-.13em o\kern-.07 em C\kern-.13em
o\kern-.15em A}}

\begin{document}

\renewcommand*{\today}{2024-03-28} 


\title{Computing the Determinant of a Dense Matrix over $\ZZ$}

\author{
  {John Abbott \qquad Claus Fieker}\\
\small{Rheinland-Pf\"alzische Technische Universit\"at Kaiserslautern}\\
\texttt{John.Abbott@rptu.de, claus.fieker@rptu.de}
}
%
%
%
\maketitle              

\begin{abstract}
  We present a new, practical algorithm for computing the determinant
  of a non-singular dense, uniform matrix over $\ZZ$; the aim is to
  achieve better practical efficiency, which is always at least as good as currently known
  methods.  The algorithm uses randomness internally, but the result
  is guaranteed correct.  The main new idea is to use a modular HNF in
  cases where the Pauderis--Storjohann HCOL method performs poorly.
  The algorithm is implemented in OSCAR~1.0.
\end{abstract}



MSC-2020: {15--04,  15A15,  15B36,  11C20}






\section{Introduction}

We present a new, practical algorithm for computing the non-zero determinant
of a matrix over $\ZZ$.  The algorithm is fully general but is intended for
dense, uniform matrices (see Section~\ref{sec:UniformDense}) We exclude
matrices with zero determinant since verifying that a matrix has zero
determinant involves different techniques; and with high probability we can
quickly detect whether the determinant is zero (\eg~checking modulo two
random 60-bit primes).

The best algorithm for computing the determinant depends on a number of
factors: \eg~dimension of the matrix, size of the entries (in relation to the dimension),
size of the determinant (usually not known in advance).
Currently, the heuristic algorithm in~\cite{PauderisStorjohann2013} is
the best choice for large dimension, dense, uniform matrices which have ``very few'', \ie~$O(1)$,
non-trivial Smith invariant factors: the ``expected''
complexity is then $O(n^3 \log(n)(\log n + \log \maxentry{A})^2)$.  The
extra condition on the Smith invariant factors very likely holds for
matrices with entries chosen from a uniform distribution over an
interval of width $\Omega(n)$~---~see~\cite{EberlyGiesbrechtVillard-2000}.  The
conclusion of~\cite{PauderisStorjohann2013} gives an indication of the
difficulty in assessing the theoretical complexity in the case that the number of
non-trivial Smith invariant factors is not $O(1)$; they do nevertheless indicate
good performance in practice when there are few non-trivial Smith invariant factors.

The main idea in~\cite{PauderisStorjohann2013} is the HCOL step which
``divides'' the working matrix by another matrix in Hermite normal form,
eventually arriving at a unimodular matrix, which is then verified to
be unimodular.  There is at least one HCOL step for each non-trivial Smith invariant factor.
We address the case where there are many non-trivial Smith invariant factors, and the
largest of these is not too large.  Instead of several rounds of division by HCOL
matrices, we do a single division by a modular HNF: the choice of strategy depends
on the size of the denominator of the solution a random linear system (as in~\cite{AbbottBronsteinMulders1999}).

\subsection{Uniform dense matrices}
\label{sec:UniformDense}

For many sorts of special matrix there are specific, efficient
algorithms for computing the determinant.  However, it can be hard to
recognize whether a given matrix belongs to one of these special
classes: \eg~deciding if a matrix is permuted triangular is NP-complete~\cite{FertinRusuVialette2015}.
We shall concern ourselves with ``uniform dense'' matrices
of integers whose \textit{entries are mostly of the same size,} 
and make no attempt to recognize any special structure.

\subsubsection*{Acknowledgements}

Both authors are supported by the Deutsche Forschungsgemeinschaft,
specifically via Project-ID~286237555~--~TRR 195

\section{Notation, Terminology, Preliminaries}

Here we introduce the notation and terminology we shall use.

\begin{definition}
  We define the \textbf{entrywise maximum} of $A \in \Mat_{r \times c}(\CC)$ as
  \[
\maxentry{A} \;= \; \max \{ A_{ij} \mid 1 \le i \le r \text{ and } 1 \le j \le c \}
\]
This is useful for specifying complexity.  We also have the following relations $\maxentry{A} \,\le\, ||A||_1 \,\le\, c \maxentry{A}$
and $\maxentry{A} \,\le\, ||A||_\infty \,\le\, r \maxentry{A}$.
\end{definition}

\begin{definition}
We say that $A \in \Mat_{n \times n}(\ZZ)$ is \textbf{unimodular}
  iff $\det(A) = \pm 1$.
\end{definition}

\begin{remark}
  In~\cite{PauderisStorjohann2012} there is an efficient algorithm to verify that
  a matrix is unimodular: it can also produce a ``certificate'' in the form of a
  ``sparse product'' representing the inverse of the matrix.
  Its complexity is $O(n^\omega \log n \, M(\log n + \log \maxentry{A}))$ where
  $\omega$ is the exponent of matrix multiplication.
\end{remark}

\begin{definition}
\label{def:Hadamard-bound}
Let $A \in \Mat_{n \times n}(\CC)$.  Then \textbf{Hadamard's row bound} for
the determinant is
\[
\Hrow(A) \;=\; \prod_{i=1}^n r_i
\]
where $r_i \in \RR_{\ge 0}$ satisfying $r^2_i = \sum_{j=1}^n |A_{ij}|^2$ is
the ``euclidean length'' of the $i$-th row of $A$.  Clearly $\Hrow(A) \le n^{n/2} \maxentry{A}^n$.
One may analogously define \textbf{Hadamard's column bound}, $\Hcol(A)$.
We may combine these to obtain
\[
H(A) \;=\; \min(\Hrow(A),\, \Hcol(A)) \;\ge\; |\det(A)|
\]
\end{definition}

\begin{remark}
  If $A \in \Mat_{n \times n}(\ZZ)$ and all entries are of similar magnitude then
  $\Hrow(A)$ and $\Hcol(A)$ are typically ``of similar size'', so there may be
  little benefit in computing both of them.  Complexity is $\softO(n^2 \log \maxentry{A})$.
\end{remark}

\begin{example}
\label{ex:Hadamard-bound}
The Hadamard row bound, $\Hrow(A)$, equals $|\det A|$ iff the rows of $A$ are
mutually orthogonal; similarly for $\Hcol(A)$.  For instance, if $A$ is
diagonal or a Hadamard matrix then both bounds are equal to the determinant.

In general, neither Hadamard bound is exact: for instance, there exist
unimodular matrices with arbitrarily large entries, such as
$\left(\begin{smallmatrix}F_{k-1}&F_k\\ F_k&F_{k+1}\end{smallmatrix}\right)$ where $F_k$ is
the $k$-th element of the Fibonacci sequence.  Also, if $n > 1$ and all entries
$A_{ij} = 1$ then $\det(A)=0$ but $H(A) = n^{n/2}$.
\end{example}

\begin{definition}
\label{def:row-HNF}
Let $A \in \Mat_{n \times m}(\ZZ)$.
The \textbf{row Hermite Normal Form} (abbr.~\textbf{row-HNF}) of $A$
is $H = U_L \, A$ where $U_L \in \Mat_{n \times n}(\ZZ)$ is
invertible, and $H \in \Mat_{n \times m}(\ZZ)$ is in ``upper
triangular'' row echelon form with positive pivots and such that for
each pivot column $j$ and each row index $i<j$ we have $0 \le H_{ij} <
H_{ii}$~---~here we adopt the usual left-to-right convention for the columns.
Observe that the row-HNF is an echelon $\ZZ$-basis for the $\ZZ$-module
generated by the rows of $A$.

The \textbf{col-HNF} may be defined analogously, but is not needed in this article.
For brevity we shall write just \textbf{HNF} to refer to the ``row'' version.
\end{definition}

\begin{definition}
\label{def:modular-HNF}  
Let $d \in \ZZ_{\neq 0}$.  Then the \textbf{$d$-modular row-HNF}
comprises the first $m$ rows of the row-HNF of
$\bigl( \begin{smallmatrix} d I_m \\ A \end{smallmatrix}\bigr)$.  This
is a full-rank, upper triangular matrix in $\Mat_{m \times m}(\ZZ)$
whose diagonal entries divide $d$.  Clearly the $d$-modular HNF is
an echelon basis for the $\ZZ$-module sum of $d\, \ZZ^m$ with the
module generated by the rows of $A$.

The \textbf{$d$-modular col-HNF} may be defined analogously, but is not needed in this article.
For brevity we shall write just \textbf{modular HNF} to refer to the
``row'' version (with the usual left-to-right convention for the
columns).
\end{definition}

\begin{remark}
The modular HNF can be computed essentially by following the
standard HNF algorithm and reducing values modulo $d$~---~except we
must not reduce the rows of $d\, I_m$.  Thus computation of the
modular HNF has bit complexity $O(n^3 (\log d)^{1+\epsilon} + n^2 \log \maxentry{A})$
assuming ``soft linear'' basic arithmetic (incl.~$\gcd$ computation).
\end{remark}

\begin{remark}
\label{rmk:hnf-diagonal}
Let $A \in \Mat_{n \times n}(\ZZ)$ be non-singular, let $p$ be a prime
dividing $\det(A)$, and let $s_n$ be the greatest Smith invariant factor of
$A$; so $s_n \neq 0$ since $A$ is non-singular.  Then for any exponent $e \in
\NN_{>0}$ every diagonal element of the $p^e$-modular HNF of $A$ divides
$p^{\min(e,k)}$ where $k = \valuation_p(s_n)$.  Also the product of the
diagonal elements divides $\det(A)$.
\end{remark}

\begin{example}
Let $d = 8$ and $M \;=\; \left(\begin{smallmatrix} 3 & -5 &  7 \\
    1 &  1 & -7 \\
    1 &  9 &  5 \\
  \end{smallmatrix}\right)$.
  Then the $d$-modular HNF is $H \;=\; \left(\begin{smallmatrix} 1 & 1 &  1 \\
    0 &  8 & 0 \\
    0 &  0 & 4 \\
  \end{smallmatrix}\right)$.
From the diagonal of $H$ we conclude that $32$ divides $\det(A)$; indeed $\det(A) = -320$.
Observe that $A  H^{-1} = \left(\begin{smallmatrix}3&-1&1\\ 1&0&-2\\ 1&1&1\end{smallmatrix}\right)$
  has integer entries~---~naturally, since the $\ZZ$-module generated by the rows of $A$ is a
  submodule of that generated by the rows of $H$.  Clearly $\det(A) / \det(H) = \det(A H^{-1}) \in \ZZ$.
\end{example}

\begin{definition}
\label{def:SNF}
Let $A \in \Mat_{n \times m}(\ZZ)$.
The \textbf{Smith Normal Form} of $A$ is $S = U_L \, A \, U_R$ where
$U_L \in \Mat_{n \times n}(\ZZ)$ and $U_R \in \Mat_{m \times m}(\ZZ)$
are invertible, and $S \in \Mat_{n \times m}(\ZZ)$ is diagonal such
that each $S_{i+1,i+1}$ is a multiple of $S_{ii}$.
By convention, the diagonal entries $S_{ii}$
are non-negative, and are known variously as \textit{elementary
  divisors, invariant factors} or \textit{Smith invariants}; we call
any $S_{ii} \neq 1$ a \textbf{non-trivial Smith invariant factor}.  We
use the standard abbreviation~\textbf{SNF}, and shall write simply
$s_j$ for the $j$-th Smith invariant factor.
\end{definition}

\begin{remark}
  With the non-negativity convention the SNF is unique; in contrast,
  the matrices $U_L$ and $U_R$ are not unique.  Be aware that some
  authors reverse the divisibility criterion of the diagonal elements;
  in which case the tuple of non-zero Smith invariants is reversed.
\end{remark}

\subsection{Unimodular matrices}
\label{sec:unimodular-matrix}

For any square $n \times n$ matrix $M$ we write $\adj(M)$ for its
adjoint (aka.~adjugate); it is well-known that $M \adj(M) = \det(M)
I_n$.  Thus for every unimodular matrix $U$ we have $U^{-1} = \pm
\adj(U)$.  The entries of $\adj(U)$ are determinants of certain
minors; this suggests that there could exist unimodular matrices whose
inverses contain large entries.  Indeed, some explicit families were
given in~\cite{NishiRumpOishi2011}: we recall one such family:  let $N \in \ZZ$ be ``large'', and let
\[
U =
\begin{pmatrix}
  1 & \epsilon & \epsilon & \cdots & \epsilon & N \\
  0 & 1 & 0 & \cdots & 0 & 0 \\
  0 & N & 1 & \cdots & 0 & 0 \\
  0 & \epsilon & N & \cdots & 0 & 0 \\
  0 & \epsilon & \epsilon & \cdots & 0 & 0 \\
  \vdots & \vdots & \vdots & \ddots & \vdots & \vdots \\
  0 & \epsilon & \epsilon & \cdots & N & 1 \\
\end{pmatrix}
\]
where taking $\epsilon < N/\sqrt{n}$ ensures that the inverse contains
at least one entry whose size is near the limit predicted by
Hadamard's bound; we may also ``perturb slightly'' the individual $N$
and $\epsilon$ values in the matrix.  Such matrices are almost worst cases
for the Pauderis--Storjohann algorithm for unimodularity verification
recalled in Section~\ref{sec:UnimodVerif}.

\subsection{Unimodularity Verification}
\label{sec:UnimodVerif}

An important step in obtaining a guaranteed result from the HCOL algorithm is verification
that a matrix is unimodular: an efficient algorithm with good
asymptotic complexity to achieve this was presented
in~\cite{PauderisStorjohann2012}.  The key to its speed is
\textit{double-plus-one lifting}~---~a variant of quadratic Hensel
lifting which cleverly limits entry growth while computing the
inverse $p$-adically.

The algorithm has worst-case performance
when the input matrix is not unimodular; so we want to avoid
applying it unless we are ``quite certain'' that the matrix is
indeed unimodular.  Also a unimodular matrix whose inverse has
large entries (\eg~the family from Section~\ref{sec:unimodular-matrix}) leads to
nearly worst-case performance.

\section{Determinant Algorithms for Integer Matrices}
\label{sec:algm-det-recall}

There are currently three practical and efficient algorithms for
computing the guaranteed determinant of an integer matrix:
multi-modular chinese remaindering (CRT), solving a random linear
system followed by some CRT steps~\cite{AbbottBronsteinMulders1999},
and the HCOL method~\cite{PauderisStorjohann2013}.  The first two
algorithms use a bound for the determinant as guarantee, whereas the
HCOL method uses \textit{unimodular verification} (see
Section~\ref{sec:UnimodVerif}).  The HCOL method has better asymptotic
complexity (and practical performance) provided that there are only
very few non-trivial Smith invariant factors.  Our new algorithm in
Section~\ref{sec:algm-det-new} addresses the case where there are
more than a few non-trivial Smith invariant factors by dividing by a
modular HNF.

\subsection{Using the modular HNF}
\label{sec:using-modular-hnf}

The cost of computing a modular HNF depends (softly linearly) on the
size of the modulus.  If the matrix $A$ has several non-trivial Smith
invariant factors then $s_n$ is likely a not-too-large factor of $\det(A)$,
and we can obtain a ``large factor'' of $s_n$ by solving a random
linear system as in~\cite{AbbottBronsteinMulders1999}~---~the factor,
$d$, appears as the common denominator of the solution.  If $d$ is ``large'',
we just do an HCOL iteration.  But if $d$ is small, we compute the
$d$-modular HNF, $H_d$, and replace $A \leftarrow A H_d^{-1} \in \Mat(\ZZ)$
noting the factor $\det(H_d)$, which is just the product of the diagonal
entries~---~this effectively condenses several HCOL iterations into
a single step. Moreover, after updating, $\det(A)$ is very likely to be
small; to make a precise probabilistic claim we would have to know
the distribution of matrices $A$ with which we compute.
Given $A$ and $d$ we can estimate quickly and reasonably accurately
how long it will take to compute $H_d$.

\section{Algorithm for Computing Determinant}
\label{sec:algm-det-new}

We present an efficient method for computing
the determinant of non-singular $A \in \Mat_{n \times n}(\ZZ)$ with $n > 2$;
we ignore any ``special structure'' $A$ may have
(and assume it was detected and handled separately during preprocessing).

\begin{itemize}
\item [(0)] Input matrix $A \in \Mat_{n \times n}(\ZZ)$ \hfill $\longleftarrow$ \textit{may be modified during algorithm}
\item [(1)] Let $e \leftarrow \lceil \log_2 \maxentry{A} \rceil$
\item [(2)] Let $h = 1+\lceil \log_2 H \rceil$ where $H \ge |\det(A)|$, \eg~Hadamard's bound.
\item [(3)] Using CRT, compute $d \leftarrow \det(A) \mod m$ where $m$ is
 a product of wordsize-bit primes; stop when either $\log_2 m > e$ or $\log_2(|d|) +60 < \log_2 m$
\item [(4)] $D \leftarrow 1$ \hfill $\longleftarrow$ \textit{always a factor of $\det(A)$}
\item [(5)]  \textbf{Main loop}
\item [(5.1)]  if $|d| = 1$ and $A$ is verified as unimodular then return $D$
\item [(5.2)]  if $60+ \log_2 |d| < \log_2 m$  and $\log_2 |d| <$ HNF threshold
\item [(5.2.1)]  Compute $|d|$-modular HNF, $H_d$; set $D_H = \det(H_d)$ and $A_{new} \leftarrow A H_d^{-1}$
\item [(5.2.2)]  If $D_H \neq 1$ update $A \leftarrow A_{new}$ and $D \leftarrow D_H D$; goto \textbf{Main loop}
\item [(5.3)] Solve linear system $A x = b$ where $b \in \Mat_{n \times 1}(\ZZ)$ is chosen randomly
\item [(5.4)] Let $D_x$ be the common denominator of the solution $x$
\item [(5.5)] if $h - \log_2(D_x d)$ less than CRT threshold
\item [(5.5.1)]  Continue chinese remaindering from step~(3) until $\log_2(m) > h/D_x$.
  \item [(5.5.2)] Return $D_x d$ where $d$ is the symmetric remainder of $\det(A) \mod m$.
\item [(5.6)] if $\log_2(D_x) <$ HNF threshold
\item [(5.6.1)]  Compute $D_x$-modular HNF, $H_{D_x}$; set $D_H = \det(H_{D_x})$
\item [(5.6.2)]  Update $A \leftarrow A H_{D_x}^{-1}$ and $D \leftarrow D_H D$; goto \textbf{Main loop}
  \item [(5.7)] From $x$ compute the ``HCOL'' matrix $H_x$
  \item [(5.8)] Update $D \leftarrow D_x D$ and $A \leftarrow A H_x^{-1}$; goto \textbf{Main loop}
\end{itemize}

\section{Practical tests}
\label{sec:timings}

We give some timings to illustrate that our method can be substantially better than
other known practical algorithms.  The examples were constructed so that the new
HNF-reduction will always take place (otherwise our method is merely a refined version of that
presented in~\cite{PauderisStorjohann2013}).  To simplify the presentation we vary
just a single parameter: the matrix dimension.  Each $n \times n$ matrix is
constructed to have $n/2$ non-trivial Smith invariant factors (being a random 11-bit
prime number), and entries of size roughly 1000 bits.

\begin{table}[!ht]
  \label{tbl:Timings}
    \centering
    \caption{Timings for determinant computation}
    \renewcommand{\arraystretch}{1.1}
    \begin{tabular}{|c||c|c|c|c|c|}
        \hline
        $n$ & New   & HCOL   & ABM  &  NTL &  CoCoA  \\
        \hline
        30  & 0.03  &  0.13  & 0.06 &  0.04 &  0.05 \\
        60  & 0.16  &  2.85  & 0.80 &  0.24 &  0.41 \\
        120 & 0.88  &  5.7   & 3.9  &  2.3  &  4.0  \\
        240 & 3.6   &  44    & 42   &  19   &  41   \\
        480 & 14.7  &  346   & 390  &  220  &  500  \\
        960 & 67    &  3570   & 4600   &  2640 &  6500  \\
        \hline
    \end{tabular}
\end{table}

Our new algorithm is already the fastest with relatively modest $30 \times 30$ matrices,
and the advantage becomes more marked as the matrix dimension increases.  The ``HCOL'' method
from~\cite{PauderisStorjohann2013} performs poorly here because our test cases were chosen
to exhibit this.  NTL~\cite{NTL} and CoCoA~\cite{CoCoA} both use chinese remaindering,
but NTL clearly has a more refined implementation.
ABM refers to the native \texttt{det} function of OSCAR which just
delegates the computation to FLINT~\cite{FLINT} which uses the method
from~\cite{AbbottBronsteinMulders1999}; for these test cases, the algorithm essentially reduces to chinese remaindering.

\section{Conclusion}

We have presented a new, practical algorithm for computing determinant of a matrix with
integer entries which exhibits good performance already for modestly sized matrices.
The new algorithm has been implemented as part of the system OSCAR~\cite{OSCAR}, and will be part
of the next major release.

\appendix


\bigskip\goodbreak

\bibliographystyle{alpha}
\bibliography{2024-03-determinant}

\goodbreak

\end{document}